\documentclass[11pt]{amsart}

\usepackage{amsmath}
\usepackage{amscd}
\usepackage{amssymb}
\usepackage{latexsym}
\usepackage{stmaryrd} 
\usepackage{mathbbol} 

\usepackage[arrow,curve,matrix,tips,2cell]{xy}
  \SelectTips{eu}{10} \UseTips
  \UseAllTwocells
\newtheorem{theorem}{Theorem}[section]
\newtheorem*{theorem*}{Theorem}
\newtheorem{lemma}[theorem]{Lemma}
\newtheorem*{lemma*}{Lemma}
\newtheorem{corollary}[theorem]{Corollary}
\newtheorem*{corollary*}{Corollary}
\newtheorem{proposition}[theorem]{Proposition}

\newtheorem{remark}[theorem]{Remark}
\newtheorem{definition}[theorem]{Definition}

\newcommand{\bgl}{\begin{equation}} 
\newcommand{\egl}{\end{equation}}
\newcommand{\bgloz}{\begin{equation*}} 
\newcommand{\egloz}{\end{equation*}}
\newcommand{\bgln}{\begin{eqnarray}} 
\newcommand{\egln}{\end{eqnarray}}
\newcommand{\bglnoz}{\begin{eqnarray*}} 
\newcommand{\eglnoz}{\end{eqnarray*}}
\newcommand{\btheo}{\begin{theorem}}
\newcommand{\etheo}{\end{theorem}}
\newcommand{\btheooz}{\begin{theorem*}}
\newcommand{\etheooz}{\end{theorem*}}
\newcommand{\blemma}{\begin{lemma}}
\newcommand{\elemma}{\end{lemma}}
\newcommand{\blemmaoz}{\begin{lemma*}}
\newcommand{\elemmaoz}{\end{lemma*}}
\newcommand{\bproof}{\begin{proof}}
\newcommand{\eproof}{\end{proof}}
\newcommand{\bbew}{\begin{beweis}}
\newcommand{\ebew}{\end{beweis}}
\newcommand{\bremark}{\begin{remark}\em}
\newcommand{\eremark}{\end{remark}}
\newcommand{\bdefin}{\begin{definition}}
\newcommand{\edefin}{\end{definition}}
\newcommand{\bprop}{\begin{proposition}}
\newcommand{\eprop}{\end{proposition}}
\newcommand{\bcor}{\begin{corollary}}
\newcommand{\ecor}{\end{corollary}}
\newcommand{\bcoroz}{\begin{corollary*}}
\newcommand{\ecoroz}{\end{corollary*}}
\newcommand{\bfa}{\begin{cases}} 
\newcommand{\efa}{\end{cases}}
\newcommand{\cI}{\mathcal I}
\newcommand{\cJ}{\mathcal J}

\newcommand{\cL}{\mathcal L}

\newcommand{\cO}{\mathcal O}
\newcommand{\cP}{\mathcal P}

\newcommand{\cR}{\mathcal R}

\def\Nz{\mathbb{N}}

\def\Zz{\mathbb{Z}}

\def\1z{\mathbb{1}}
\newcommand{\fA}{\mathfrak A}

\newcommand{\fI}{\mathfrak I}

\newcommand{\mfa}{\mathfrak a}

\newcommand{\mfk}{\mathfrak k}
\newcommand{\mfp}{\mathfrak p}
\newcommand{\mfq}{\mathfrak q}
\newcommand{\an}[1]{``#1''} 

\newcommand{\ma}{\mapsto} 
\newcommand{\mafr}{\mapsfrom} 
\newcommand\into{\hookrightarrow} 
\newcommand{\Rarr}{\Rightarrow} 
\newcommand{\LRarr}{\Leftrightarrow} 

\def\SEMI{\mbox{$\times\kern-2pt\vrule height5pt width.6pt \kern3pt $}}

\newcommand{\Spec}{{\rm Spec\,}} 
\newcommand{\Prim}{{\rm Prim\,}} 

\renewcommand{\ker}{{\rm ker}\,}

\newcommand{\reg}{^\times} 
\newcommand{\abs}[1]{\lvert#1\rvert} 
\newcommand{\defeq}{\mathrel{:=}} 

\newcommand{\dop}{\text{: }} 
\newcommand{\falls}{\text{ if }} 
\newcommand{\sonst}{\text{ else}} 
\newcommand{\fa}{\text{ for all }} 
\newcommand{\ilim}{\varinjlim} 
\newcommand{\plim}{\varprojlim} 

\newcommand{\height}{\operatorname{ht}} 
\newcommand{\tfin}{\text{fin}} 
\newcommand{\tinf}{\text{inf}}
\newcommand{\supp}{\operatorname{supp}} 

\newcommand{\lge}{\left\{} 
\newcommand{\rge}{\right\}} 
\newcommand{\lru}{\left(} 
\newcommand{\rru}{\right)} 
\newcommand{\lsp}{\left\langle} 
\newcommand{\rsp}{\right\rangle} 
\newcommand{\rukl}[1]{\lru #1 \rru} 
\newcommand{\gekl}[1]{\lge #1 \rge} 
\newcommand{\spkl}[1]{\lsp #1 \rsp} 
\newcommand{\menge}[2]{\gekl{ #1 \dop #2 }} 
\newcommand{\ping}{P \subseteq G}

\begin{document}

\title{Semigroup C*-algebras of $ax+b$-semigroups}

\author{Xin Li}

\address{Xin Li, Department of Mathematics, Westf{\"a}lische Wilhelms-Universit{\"a}t M{\"u}nster, Einsteinstra{\ss}e 62, 48149 M{\"u}nster, Germany}
\email{xinli.math@uni-muenster.de}


\thanks{\scriptsize{Research supported by the ERC through AdG 267079.}}

\begin{abstract}
We study semigroup C*-algebras of $ax+b$-semigroups over integral domains. The goal is to generalize several results about C*-algebras of $ax+b$-semigroups over rings of algebraic integers. We prove results concerning K-theory and structural properties like the ideal structure or pure infiniteness. Our methods allow us to treat $ax+b$-semigroups over a large class of integral domains containing all noetherian, integrally closed domains and coordinate rings of affine varieties over infinite fields.
\end{abstract}

\maketitle


\setlength{\parindent}{0pt} \setlength{\parskip}{0.5cm}

\section{Introduction}

Given an integral domain, let us form the $ax+b$-semigroup over this ring and consider the C*-algebra generated by the left regular representation of this semigroup. The particular case where the integral domain is given by the ring of algebraic integers in a number field has been investigated in \cite{C-D-L}, \cite{C-E-L1}, \cite{C-E-L2}, \cite{E-L}, \cite{Li4} and also \cite{La-Ne}, \cite{L-R}. The C*-algebras of such $ax+b$-semigroups turn out to have intriguing structural properties as well as very interesting K-theory. In a bigger context, the study of such C*-algebras has initiated the investigation of much more general semigroup C*-algebras and has led to a better understanding of their structure (see \cite{Li2}, \cite{Li3}, \cite{C-E-L1}, \cite{C-E-L2}, \cite{Nor}).

In this paper, our goal is to generalize several results about the C*-algebras of these $ax+b$-semigroups from rings of algebraic integers (\cite{C-E-L1}, \cite{C-E-L2}, \cite{E-L}) to much more general classes of integral domains. This generalization process reveals which properties of the rings are responsible for which properties of the semigroup C*-algebras. Therefore, our present work provides a better understanding of the original results in the case of rings of algebraic integers. At the same time, our methods and results considerably enlarge the source of (tractable) examples of semigroup C*-algebras which come from rings.

Let us now formulate our main results. Given an integral domain $R$, we consider both the left semigroup C*-algebra $C^*_{\lambda}(R \rtimes R\reg)$ generated by the left regular representation of the $ax+b$-semigroup $R \rtimes R\reg$, as well as the right semigroup C*-algebra $C^*_{\rho}(R \rtimes R\reg)$ generated by the right regular representation of $R \rtimes R\reg$. A main theme is to compare $C^*_{\lambda}(R \rtimes R\reg)$ and $C^*_{\rho}(R \rtimes R\reg)$.

Our first main result is concerned with K-theory:
\btheo
\label{main-K}
Let $R$ be a countable Krull ring with group of multiplicative units $R^*$ and divisor class group $C(R)$. For every $\mfk \in C(R)$, let $\mfa_\mfk$ be a divisorial ideal which represents $\mfk$. Then
$$
  K_*(C^*_{\lambda}(R \rtimes R\reg)) \cong \bigoplus_{\mfk \in C(R)} K_*(C^*(\mfa_\mfk \rtimes R^*)),
$$
$$
  K_*(C^*_{\rho}(R \rtimes R\reg)) \cong \bigoplus_{\mfk \in C(R)} K_*(C^*(\mfa_{\mfk^{-1}} \rtimes R^*)).
$$
In particular, $K_*(C^*_{\lambda}(R \rtimes R\reg)) \cong K_*(C^*_{\rho}(R \rtimes R\reg))$.
\etheo
This first result generalizes the K-theoretic computations in \cite[\S~8]{C-E-L1} and \cite[\S~6.4]{C-E-L2}. The notion of a Krull ring is explained in \S~\ref{sec-Krull}. Every noetherian, integrally closed integral domain is a Krull ring. For instance, for every affine scheme which is integral, noetherian and regular, the corresponding coordinate ring is a Krull ring. For such schemes, the divisor class group can be identified with the group of Weil or Cartier divisors (these notions coincide for such regular schemes).

Secondly, we are able to compute the primitive ideal space of $C^*_{\lambda}(R \rtimes R\reg)$. Let $R$ be a Krull ring. Let $\cP(R)$ be the set of prime ideals of $R$ which are of height $1$, and let $\cP_{\tfin}(R) \defeq \menge{\mfp \in \cP(R)}{[R:\mfp^{(i)}] < \infty \fa i \in \Nz}$ and $\cP_{\tinf}(R) \defeq \cP(R) \setminus \cP_{\tfin}(R)$. Here $\mfp^{(i)}$ is the $i$-fold product of $\mfp$ with itself in the monoid of divisorial ideals (see \S~\ref{sec-Krull}).
\btheo
\label{main-primitive}
Assume that $R$ is a countable Krull ring and that $\cP(R)_{\tinf} \neq \emptyset$ or that $\cP(R)_{\tfin}$ is infinite. Then there is an order-preserving homeomorphism $\Prim(C^*_{\lambda}(R \rtimes R\reg)) \cong 2^{\cP_{\tfin}(R)}$, where $2^{\cP_{\tfin}(R)}$ is the power set of $\cP_{\tfin}(R)$ equiped with the power-cofinite topology. All the orders are given by inclusion (of primitive ideals or subsets).
\etheo
This generalizes \cite[Theorem~3.6]{E-L}.

Thirdly, we show that for certain rings $R$, $C^*_{\lambda}(R \rtimes R\reg)$ is strongly purely infinite, i.e., $\cO_{\infty}$-absorbing. This property is of importance for the classification of C*-algebras.
\btheo
\label{main-pi}
Let $R$ be a countable integral domain. Suppose that $R$ is not a field and that the Jacobson radical of $R$ (the intersection of all maximal ideals of $R$) vanishes. Then $C^*_{\lambda}(R \rtimes R\reg)$ is purely infinite and has the ideal property, and is therefore strongly purely infinite ($C^*_{\lambda}(R \rtimes R\reg) \cong C^*_{\lambda}(R \rtimes R\reg) \otimes \cO_{\infty}$).
\etheo
This result generalizes \cite[Theorem~8.2.4]{C-E-L1}.

Note that while in K-theory, there is no difference between left and right semigroup C*-algebras for $ax+b$-semigroups over Krull rings, the analogues of Theorem~\ref{main-primitive} and Theorem~\ref{main-pi} cannot be true for right semigroup C*-algebras of $ax+b$-semigroups.

Apart from that, we also have the following result:
\btheo
\label{main-inffield}
Let $R$ be a countable integral domain which contains an infinite field. Then $C^*_{\lambda}(R \rtimes R\reg)$ is a unital UCT Kirchberg algebra.

To describe the K-theory of $C^*_{\lambda}(R \rtimes R\reg)$ for such rings, let $Q$ be the quotient field of $R$, set $\cI(R \subseteq Q) \defeq \menge{(x_1 \cdot R) \cap \dotso (x_n \cdot R)}{x_i \in Q\reg}$ and for $I \in \cI(R \subseteq Q)$, let $Q\reg_I \defeq \menge{a \in Q\reg}{aI = I}$. Let $Q\reg \backslash \cI(R \subseteq Q)$ be the set of orbits for the canonical multiplicative action of $Q\reg$ on $\cI(R \subseteq Q)$. Then
$$
  K_*(C^*_{\lambda}(R \rtimes R\reg)) \cong \bigoplus_{[I] \in Q\reg \backslash \cI(R \subseteq Q)} K_*(C^*_{\lambda}(I \rtimes Q\reg_I)).
$$
The $K_0$-class of the unit of $C^*_{\lambda}(R \rtimes R\reg)$ corresponds to the $K_0$-class of the unit of $C^*_{\lambda}(R \rtimes R^*)$ for the orbit $[R] \in Q\reg \backslash \cI(R \subseteq Q)$, i.e., the orbit of the prinipal fractional ideals.
\etheo
This result applies to all coordinate rings of affine varieties over infinite fields. Also, note that by \cite[Chapter~8]{Ror}, unital UCT Kirchberg algebras are completely classified by their K-theory together with the position of the $K_0$-class of the unit.

Let us briefly explain what sort of methods we use to obtain our main results: To prove Theorem~\ref{main-K} and Theorem~\ref{main-inffield}, we study the independence condition (see \S~\ref{sec-intdom}) for $ax+b$-semigroups over integral domains. Once the independence condition is understood, we only have to apply results from \cite{Li1}, \cite{C-E-L1} and \cite{C-E-L2}. For Theorem~\ref{main-primitive}, the most important ingredient is a careful analysis of the spectrum of a canonical commutative subalgebra (see \S~\ref{semiC}) of our semigroup C*-algebra and its dilated version. This spectrum can be thought of as a substitute for the finite adele space in the case of rings of algebraic integers (in number fields or global fields). But it is really a modified version of the finite adele space because we systematically adjoin \an{points at infinity} in order to make sure that we get a locally compact space. This \an{local compactification procedure} might be interesting in other contexts as well. And finally, to prove Theorem~\ref{main-pi}, we basically refine the methods from \cite{C-E-L1} in order to treat more general rings than rings of algebraic integers. This refinement naturally leads to our condition involving the Jacobson radical.

This paper is structured as follows: First, we recall the construction of semigroup C*-algebras (\S~\ref{semiC}). We then consider the special case of $ax+b$-semigroups and study the independence condition (\S~\ref{sec-intdom}). Furthermore, we explain the notion of Krull rings in \S~\ref{sec-Krull} and prove that independence holds in the case of Krull rings. In \S~\ref{functor}, we briefly discuss functorial properties of C*-algebras attached to $ax+b$-semigroups. We determine K-theory and prove Theorem~\ref{main-K} in \S~\ref{sec-K}. In \S~\ref{sec-prim}, we study the ideal structure and prove Theorem~\ref{main-primitive}. And pure infiniteness is discussed in \S~\ref{sec-pi} where we prove Theorem~\ref{main-pi} and Theorem~\ref{main-inffield}. Finally, we give an example of an integral domain which does not satisfy independence in \S~\ref{ind-fails}.

{\bf Acknowledgement:} I would like to thank Torsten Schoeneberg for helpful comments about Krull rings.

\section{Preliminaries}

\subsection{Semigroup C*-algebras}
\label{semiC}

Let us recall the construction of semigroup C*-algebras. Given a left cancellative semigroup $P$, consider the Hilbert space $\ell^2 P$ with its canonical orthonormal basis $\menge{\varepsilon_x}{x \in P}$, and define for every $p \in P$ an isometry $\lambda(p)$ by setting $\lambda(p) \varepsilon_x = \varepsilon_{px}$. As for reduced group C*-algebras, we simply take the C*-algebra generated by the left regular representation of our semigroup:
$$C^*_{\lambda}(P) \defeq C^* \rukl{\menge{\lambda(p)}{p \in P}} \subseteq \cL(\ell^2 P).$$
Similarly, given a right cancellative semigroup $P$, we can form the operators $\rho(p)$ on $\ell^2 P$ defined by $\rho(p) \varepsilon_x = \varepsilon_{xp}$, and consider the C*-algebra generated by the right regular representation of $P$:
$$C^*_{\rho}(P) \defeq C^* \rukl{\menge{\rho(p)}{p \in P}} \subseteq \cL(\ell^2 P).$$
Results about left semigroup C*-algebras can be transfered into results about right semigroup C*-algebras (and vice versa) using the relation $C^*_{\rho}(P) \cong C^*_{\lambda}(P^{\text{op}})$ for a right cancellative semigroup $P$. Here $P^{\text{op}}$ is the opposite semigroup where the multiplication is flipped.

In the analysis of these semigroup C*-algebras, two conditions play a prominent role. The first one is a condition on the ideal structure of the semigroup.
\bdefin
For a left cancellative semigroup $P$, we let
$$\cJ_{\lambda}(P) \defeq \menge{q_1^{-1} p_1 \dotsm q_n^{-1} p_n P}{p_i, q_i \in P} \cup \gekl{\emptyset}$$
be the family of constructible right ideals of $P$.
\edefin
Here, for $q \in P$ and a subset $X$ of $P$, $q^{-1} X = \menge{y \in P}{qy \in X}$.
\bdefin
We call $\cJ_{\lambda}(P)$ independent if for every $X, X_1, \dotsc, X_n \in \cJ_{\lambda}(P)$, the following holds: Whenever $X = \bigcup_{i=1}^n X_i$, then we must have $X = X_i$ for some $1 \leq i \leq n$.

We say that $P$ satisfies the left independence condition if $\cJ_{\lambda}(P)$ is independent.
\edefin
Of course, all these definitions have their right analogues:
\bdefin
Given a right cancellative semigroup $P$, let
$$\cJ_{\rho}(P) \defeq \menge{P p_n q_n^{-1} \dotsm p_1 q_1^{-1}}{p_i, q_i \in P} \cup \gekl{\emptyset}$$
be the family of constructible left ideals of $P$.
\edefin
Here, for $q \in P$ and a subset $X$ of $P$, $X q^{-1} = \menge{y \in P}{yq \in X}$.
\bdefin
We call $\cJ_{\rho}(P)$ independent if for every $X, X_1, \dotsc, X_n \in \cJ_{\rho}(P)$, the following holds: Whenever $X = \bigcup_{i=1}^n X_i$, then we must have $X = X_i$ for some $1 \leq i \leq n$.

We say that $P$ satisfies the right independence condition if $\cJ_{\rho}(P)$ is independent.
\edefin

The second condition is a condition on an embedding of our semigroup into a group. Assume that $P \subseteq G$ is an embedding of $P$ as a subsemigroup into a group $G$. Let $E_P$ be the orthogonal projection in $\cL(\ell^2 G)$ onto the subspace $\ell^2 P \subseteq \ell^2 G$. Moreover, let $\lambda^G$ and $\rho^G$ be the left and right regular representations of $G$ on $\ell^2 G$.
\bdefin
We say that $\ping$ satisfies the left Toeplitz condition (or simply that $\ping$ is left Toeplitz) if for every $g \in G$ with $E_P \lambda^G_g E_P \neq 0$, there exist $p_1, q_1, \dotsc, p_n, q_n$ in $P$ such that $E_P \lambda^G_g E_P = \lambda(p_1)^* \lambda(q_1) \dotsm \lambda(p_n)^* \lambda(q_n)$.
\edefin
\bdefin
We say that $\ping$ satisfies the right Toeplitz condition (or simply that $\ping$ is right Toeplitz) if for every $g \in G$ with $E_P \rho^G_g E_P \neq 0$, there exist $p_1, q_1, \dotsc, p_n, q_n$ in $P$ such that $E_P \rho^G_g E_P = \rho(p_1)^* \rho(q_1) \dotsm \rho(p_n)^* \rho(q_n)$.
\edefin
If $P$ is subsemigroup of a group $G$, the following notations will be helpful:
\bdefin
The family of constructible right $P$-ideals in $G$ is given by
$$\cJ_{\lambda}(\ping) \defeq \menge{\bigcap_{i=1}^n g_i \cdot P}{g_1, \dotsc, g_n \in G} \cup \gekl{\emptyset}.$$
\edefin
If $\ping$ is left Toeplitz, then $\cJ_{\lambda}(\ping) = \menge{g \cdot X}{g \in G, \, X \in \cJ_{\lambda}(P)}$.

We also have a right analogue:
\bdefin
The family of constructible left $P$-ideals in $G$ is given by
$$\cJ_{\rho}(\ping) \defeq \menge{\bigcap_{i=1}^n P \cdot g_i}{g_1, \dotsc, g_n \in G} \cup \gekl{\emptyset}.$$
\edefin
Again, if $\ping$ is right Toeplitz, then $\cJ_{\rho}(\ping) = \menge{X \cdot g}{g \in G, \, X \in \cJ_{\rho}(P)}$.

Note that for an element $q \in P \subseteq G$ and a subset $X$ of $P \subseteq G$, we write $q^{-1}X$ for the set $\menge{y \in P}{qy \in X}$, whereas $q^{-1} \cdot X = \menge{q^{-1}x}{x \in X}$. In general, these sets do not coincide, but they are related as follows: $q^{-1}X = (q^{-1} \cdot X) \cap P$. An analogous comment applies to the right versions of these notations.

Let us explain why the Toeplitz condition is so useful. We restrict ourselves to the left Toeplitz condition, but we have analogous results for the right version. Let $P$ be a subsemigroup of a group $G$. The semigroup C*-algebra $C^*_{\lambda}(P)$ contains a canonical commutative sub-C*-algebra $D_{\lambda}(P)$. It is given by $D_{\lambda}(P) = C^*(\menge{E_X}{X \in \cJ_{\lambda}(P)})$. Here $E_X \in \cL(\ell^2 P)$ is the orthogonal projection onto the subspace $\ell^2 X$ of $\ell^2 P$ for $X \subseteq P$. Moreover, let $D_{\lambda, \ping}$ be the smallest $G$-invariant sub-C*-algebra of $\ell^{\infty} G$ which contains $D_{\lambda}(P)$. In \cite{Li3}, $D_{\lambda, \ping}$ is denoted by $D_P^G$. $C^*_{\lambda}(P)$ sits in a canonical way in $D_{\lambda, \ping} \rtimes_r G$ (see \cite[\S~3]{Li3}), and we will from now on identify $C^*_{\lambda}(P)$ with the corresponding sub-C*-algebra of $D_{\lambda, \ping} \rtimes_r G$. The point now is that if $\ping$ is left Toeplitz, then $C^*_{\lambda}(P)$ is not only a canonical sub-C*-algebra of $D_{\lambda, \ping} \rtimes_r G$, but actually a full corner. This means that $C^*_{\lambda}(P)$ and $D_{\lambda, \ping} \rtimes_r G$ are Morita equivalent, so that (among other things) their K-theories and primitive ideal spaces coincide. Since crossed products of the from $D_{\lambda, \ping} \rtimes_r G$ have already been intensively studied, this allows us to establish many structural results about $C^*_{\lambda}(P)$. This is why the Toeplitz condition is so useful.

\subsection{The case of $ax+b$-semigroups over integral domains}
\label{sec-intdom}

We now specialize to the following situation: Let $R$ be an integral domain, by which we mean a commutative ring with unit ($0 \neq 1$) which does not have zero-divisors. Moreover, we assume that $R$ is countable so that our C*-algebras will be separable. We are interested in the left and right semigroup C*-algebras of the $ax+b$-semigroup $R \rtimes R\reg$ over $R$. By definition, $R \rtimes R\reg$ is the semidirect product of the additive group $R$ by the multiplicative semigroup $R\reg \defeq R \setminus \gekl{0}$ with respect to the multiplicative action of $R\reg$ on $R$. To be explicit, as a set, $R \rtimes R\reg$ is the direct product $R \times R\reg$, and the multiplication is given by $(b,a)(x,w) = (b+ax,aw)$ for $(b,a)$, $(x,w)$ in $R \times R\reg$.
\bdefin
The constructible ring-theoretic ideals of $R$ are given by
$$
  \cI(R) \defeq \menge{c^{-1} \rukl{\bigcap_{i=1}^n a_i R}}{a_1, \dotsc, a_n, c \in R\reg}.
$$
\edefin
Here, for $c \in R\reg$ and an ideal $I$ of $R$, we set $c^{-1}I \defeq \menge{r \in R}{cr \in I}$.
\bdefin
We call $\cI(R)$ independent if for every $I, I_1, \dotsc, I_n \in \cI(R)$, the following holds: Whenever $I = \bigcup_{i=1}^n I_i$, then we must have $I = I_i$ for some $1 \leq i \leq n$.

We say that $R$ satisfies the independence condition if $\cI(R)$ is independent.
\edefin

A simple computation shows the following
\blemma
We have $\cJ_{\lambda}(R \rtimes R\reg) = \menge{(r+I) \times I\reg}{r \in R, I \in \cI(R)} \cup \gekl{\emptyset}$ and $\cJ_{\rho}(R \rtimes R\reg) = \menge{R \times I\reg}{I \in \cI(R)} \cup \gekl{\emptyset}$, where $I\reg \defeq I \setminus \gekl{0}$.
\elemma

Let us make the following observation about the relationship between the left and right independence condition for $ax+b$-semigroups:
\blemma
\label{independence}
The following are equivalent:
\begin{itemize}
\item $\cJ_{\lambda}(R \rtimes R\reg)$ is independent,
\item $\cJ_{\rho}(R \rtimes R\reg)$ is independent,
\item $\cI(R)$ is independent.
\end{itemize}
\elemma
\bproof
It is obvious that the second and third items are equivalent. So it remains to show that $\cJ_{\lambda}(R \rtimes R\reg)$ is independent if and only if $\cI(R)$ is independent. If $\cJ_{\lambda}(R \rtimes R\reg)$ is not independent, then we have a non-trivial equation of the form
$$
  (r+I) \rtimes I\reg = \bigcup_{i=1}^n (r_i + I_i) \times I_i\reg \text{ with } (r_i + I_i) \times I_i\reg \subsetneq (r+I) \rtimes I\reg.
$$
$(r_i + I_i) \times I_i\reg \subsetneq (r+I) \rtimes I\reg$ implies that $I_i \subsetneq I$, and projecting onto the second coordinate of $R \times R\reg$, we obtain $I\reg = \bigcup_{i=1}^n I_i\reg$, hence $I = \bigcup_{i=1}^n I_i$. This means that $\cI(R)$ is not independent. Conversely, assume that $\cI(R)$ is not independent, so that we have a non-trivial equation of the form $I = \bigcup_{i=1}^n I_i$ with $I_i \subsetneq I$. By \cite[Theorem~18]{Gott}, we may assume without loss of generality that $[I:I_i] < \infty$ for all $1 \leq i \leq n$. But then we have
$$
  I \times I\reg = \bigcup_{i=1}^n \bigcup_{r + I_i \in I / I_i} (r+I_i) \times I_i\reg.
$$
This shows that $\cJ_{\lambda}(R \rtimes R\reg)$ is not independent.
\eproof

Of course, a natural question that arises at this point is: Which integral domains satisfy the independence condition? The complete answer to this question is not known to the author. However, we can give a few sufficient conditions:
\begin{itemize}
\item If for every $I$, $J$ in $\cI(R)$ with $I \subsetneq J$, we have $[J:I] = \infty$, then $\cI(R)$ is independent. This follows from \cite[Lemma~17]{Gott}. In particular, this is the case if $R$ contains an infinite field $k$, because then $J/I$ is a $k$-vector space, hence infinite.
\item If every ideal in $\cI(R)$ is an invertible ideal, then $\cI(R)$ is independent. This follows from \cite[Theorem~1.5]{B-Q}.
\item If $R$ is a Krull ring, then $\cI(R)$ is independent. This will be proven in the next paragraph.
\end{itemize}

Next, we turn to the Toeplitz condition. Let $Q$ be the quotient field of $R$. It is clear that $R \rtimes R\reg$ is in a canonical way a subsemigroup of the $ax+b$-group $Q \rtimes Q\reg$.
\blemma
\label{Toeplitz}
$R \rtimes R\reg \subseteq Q \rtimes Q\reg$ satisfies both the left and right Toeplitz conditions.
\elemma
\bproof
Since $R \rtimes R\reg$ is a left Ore semigroup with enveloping group $Q \rtimes Q\reg$, $R \rtimes R\reg \subseteq Q \rtimes Q\reg$ satisfies the left Toeplitz condition by \cite[\S~8.3]{Li3}.

To see that $R \rtimes R\reg \subseteq Q \rtimes Q\reg$ satisfies the right Toeplitz condition, we apply \cite[Proposition~6.1.7]{C-E-L2} to $\bar H = (Q,+)$, $\bar P = (Q\reg,\cdot)$, $H = (R,+)$ and $P = (R\reg,\cdot)$: For $\bar h = q \in Q\reg$, let $z \defeq q$, $h \defeq 1$. Then $R\reg \cap R\reg z^{-1} = R\reg \cap R\reg q^{-1} = \menge{x \in R\reg}{xq \in R\reg} = \menge{x \in R\reg}{xq \in R}$ as $q \neq 0$. For $\bar h = 0$, set $z \defeq 1$ and $h \defeq 0$. Then $R\reg \cap R\reg z^{-1} = R\reg = \menge{x \in R\reg}{x \cdot 0 \in R}$. Thus (6.2) from \cite[Proposition~6.1.7]{C-E-L2} holds, and we are done.
\eproof

Moreover, let us explicitly write down the constructible $R \rtimes R\reg$-ideals in $Q \rtimes Q\reg$. Let
\bgl
\label{const-ring-ideals}
  \cI(R \subseteq Q) \defeq \menge{(x_1 \cdot R) \cap \dotso (x_n \cdot R)}{x_i \in Q\reg}.
\egl
Note that for $c \in R\reg$ and $X \subseteq R$, we set $c^{-1} X = \menge{r \in R}{cr \in X}$, but ${c^{-1} \cdot X} = \menge{c^{-1} x}{x \in X}$. Moreover, note that $\cI(R) = \menge{J \cap R}{J \in \cI(R \subseteq Q)}$. Given $q \in Q$ and $J \in \cI(R \subseteq Q)$, we set $Y(q,J) \defeq \menge{(b,a) \in Q \rtimes Q\reg}{a \in J\reg, \, b \in aq + R}$. Straightforward computations show the following
\blemma
$\cJ_{\lambda}(R \rtimes R\reg \subseteq Q \rtimes Q\reg) = \menge{(q+J) \times J\reg}{q \in Q, \, J \in \cI(R \subseteq Q)} \cup \gekl{\emptyset}$ and $\cJ_{\rho}(R \rtimes R\reg \subseteq Q \rtimes Q\reg) = \menge{Y(q,J)}{q \in Q, \, J \in \cI(R \subseteq Q)} \cup \gekl{\emptyset}$.
\elemma

\subsection{Krull rings}
\label{sec-Krull}

By construction, the family $\cI(R)$ consists of integral divisorial ideals of $R$, and $\cI(R \subseteq Q)$ consists of divisorial ideals of $R$. By definition, a divisorial ideal of an integral domain $R$ is a fractional ideal $I$ that satisfies $I = (R:(R:I))$, where $(R:J) = \menge{q \in Q}{qJ \subseteq R}$. Equivalently, divisorial ideals are non-zero intersections of some non-empty family of principal fractional ideals (ideals of the form $qR$, $q \in Q$). Let $D(R)$ be the set of divisorial ideals of $R$. In our situation, we only consider finite intersections of principal fractional ideals (see \eqref{const-ring-ideals}). So in general, our family $\cI(R \subseteq Q)$ will only be a proper subset of $D(R)$.

However, for certain rings, the set $\cI(R \subseteq Q)$ coincides with $D(R)$. For instance, this happens for noetherian rings. It also happens for Krull rings. The latter have a number of additional favourable properties which are very helpful for our purposes:
\bdefin
\label{Krull}
An integral domain $R$ is called a Krull ring if there exists a family of discrete valuations $(v_i)_{i \in I}$ of the quotient field $Q$ of $R$ such that
\begin{enumerate}
\item[(K1)] $R = \menge{x \in Q}{v_i(x) \geq 0 \fa i \in I}$,
\item[(K2)] for every $0 \neq x \in Q$, there are only finitely many valuations in $(v_i)_i$ such that $v_i(x) \neq 0$.
\end{enumerate}
\edefin

The following result gives us many examples of Krull rings.
\btheo{\cite[Chapitre~VII, \S~1.3, Corollaire]{Bour2}}
A noetherian integral domain is a Krull ring if and only if it is integrally closed.
\etheo

Let us collect some basic properties of Krull rings:

\cite[Chapitre~VII, \S~1.5, Corollaire~2]{Bour2} yields
\blemma
For a Krull ring $R$, $\cI(R \subseteq Q) = D(R)$ and $\cI(R)$ is the set of integral divisorial ideals.
\elemma
Moreover, the prime ideals of height $1$ play a distinguished role in a Krull ring.
\btheo{\cite[Chapitre~VII, \S~1.6, Th\'{e}or\`{e}me~3 and Chapitre~VII, \S~1.7, Th\'{e}or\`{e}me~4]{Bour2}}
Let $R$ be a Krull ring. Every prime ideal of height $1$ of $R$ is a divisorial ideal. Let $\cP(R) = \menge{\mfp \triangleleft R \text{ prime}}{\height(\mfp)=1}$. For every $\mfp \in \cP(R)$, the localization $R_{\mfp} = (R \setminus \mfp)^{-1} R$ is a principal valuation ring. Let $v_{\mfp}$ be the corresponding (discrete) valuations of the quotient field $Q$ of $R$. Then the family $(v_{\mfp})_{\mfp \in \cP(R)}$ satisfies the conditions (K1) and (K2) from Definition~\ref{Krull}.
\etheo
\bprop{\cite[Chapitre~VII, \S~1.5, Proposition~9]{Bour2}}
\label{approx}
Let $R$ be a Krull ring and $(v_{\mfp})_{\mfp \in \cP(R)}$ be the valuations from the previous theorem. Given finitely many integers $n_1$, ..., $n_r$ and finitely many prime ideals $\mfp_1$, ..., $\mfp_r$ in $\cP(R)$, there exists $x$ in the quotient field $Q$ of $R$ with $v_{\mfp_i}(x) = n_i$ for all $1 \leq i \leq r$ and $v_{\mfp}(x) \geq 0$ for all $\mfp \in \cP(R) \setminus \gekl{\mfp_1, \dotsc, \mfp_r}$.
\eprop
Moreover, given a fractional ideal $I$ of $R$, we let $I^\sim \defeq (R:(R:I))$ be the divisorial closure of $I$. $I^\sim$ is the smallest divisorial ideal of $R$ which contains $I$. We can now define the product of two divisorial ideals $I_1$ and $I_2$ to be the divisorial closure of the (usual ideal-theoretic) product of $I_1$ and $I_2$, i.e., $I_1 \bullet I_2 \defeq (I_1 \cdot I_2)^\sim$. $D(R)$ becomes a commutative monoid with this multiplication.
\btheo{\cite[Chapitre~VII, \S~1.2, Th\'{e}or\`{e}me~1; Chapitre~VII, \S~1.3, Th\'{e}or\`{e}me~2 and Chapitre~VII, \S~1.6, Th\'{e}or\`{e}me~3]{Bour2}}
For a Krull ring $R$, $(D(R),\bullet)$ is a group, namely the free abelian group with free generators given by $\cP(R)$, the set of prime ideals of $R$ which have height $1$.
\etheo
This means that every $I \in \cI(R \subseteq Q)$ ($Q$ is the quotient field of the Krull ring $R$) is of the form $I = \mfp_1^{(n_1)} \bullet \dotsm \bullet \mfp_r^{(n_r)}$, with $n_i \in \Zz$. Here for $\mfp \in \cP(R)$ and $n \in \Nz$, we write $\mfp^{(n)}$ for $\underbrace{\mfp \bullet \dotsm \bullet \mfp}_{n \text{ times}}$, and $\mfp^{(-n)}$ for $\underbrace{\mfp^{-1} \bullet \dotsm \bullet \mfp^{-1}}_{n \text{ times}}$, where $\mfp^{-1} = (R:\mfp)$. We set for $\mfp \in \cP(R)$:
$
  v_{\mfp}(I) \defeq
  \bfa
    n_i \falls \mfp = \mfp_i, \\
    0 \falls \mfp \notin \gekl{\mfp_1, \dotsc, \mfp_r}.
  \efa
$
With this notation, we have $I = \prod_{\mfp \in \cP(R)} \mfp^{(v_{\mfp}(I))}$, where the product is taken in $D(R)$. In addition, we have for $I \in \cI(R \subseteq Q)$: $I \in \cI(R)$ $\LRarr$ $v_{\mfp}(I) \geq 0$ for all $\mfp \in \cP(R)$. And combining the last statement in \cite[Chapitre~VII, \S~1.3, Th\'{e}or\`{e}me~2]{Bour2} with \cite[Chapitre~VII, \S~1.4, Proposition~5]{Bour2}, we obtain for every $I \in \cI(R \subseteq Q)$:
\bgl
\label{I=}
  I = \menge{x \in Q}{v_{\mfp}(x) \geq v_{\mfp}(I) \text{ for all } \mfp \in \cP(R)}.
\egl
Finally, the principal fractional ideals $F(R)$ form a subgroup of $(D(R),\bullet)$ which is isomorphic to $Q\reg$. Suppose that $R$ is a Krull ring. Then the quotient group $C(R) \defeq D(R) / F(R)$ is called the divisor class group of $R$.

These were basic properties of Krull rings. We refer the interested reader to \cite[Chapitre~VII]{Bour2} or \cite{Fos} for more information.

For us, the following consequence is of particular importance:
\blemma
\label{Krull-independence}
A Krull ring satisfies the independence condition.
\elemma
\bproof
Let $R$ be a Krull ring with quotient field $Q$, and let $I$, $I_1$, ..., $I_n$ be ideals in $\cI(R)$ with $I_i \subsetneq I$ for all $1 \leq i \leq n$. Then for every $1 \leq i \leq n$, there exists $\mfp_i \in \cP(R)$ with $v_{\mfp_i}(I_i) > v_{\mfp_i}(I)$. By Proposition~\ref{approx}, there exists $x \in Q\reg$ with $v_{\mfp_i}(x) = v_{\mfp_i}(I)$ for all $1 \leq i \leq n$ and $v_{\mfp}(x) \geq v_{\mfp}(I)$ for all $\mfp \in \cP(R) \setminus \gekl{\mfp_1, \dotsc, \mfp_r}$. Thus $x$ lies in $I$, but does not lie in $I_i$ for any $1 \leq i \leq n$. Therefore, $\bigcup_{i=1}^n I_i \subsetneq I$. This shows that $\cI(R)$ is independent, as claimed.
\eproof

\section{Functorial properties}
\label{functor}

Let us discuss functorial properties of semigroup C*-algebras of $ax+b$-semigroups.
\bprop
For rings satisfying the independence condition, our assignment $R \ma C^*_{\lambda}(R \rtimes R\reg)$ is functorial for faithfully flat ring monomorphims, i.e., every faithfully flat ring monomorphism $\phi: R \into S$ induces a homomorphism $C^*_{\lambda}(R \rtimes R\reg) \to C^*_{\lambda}(S \rtimes S\reg)$ determined by $\lambda(b,a) \ma \lambda(\phi(b),\phi(a))$.

Moreover, again for rings satisfying the independence condition, our assignment $R \ma C^*_{\rho}(R \rtimes R\reg)$ is functorial for flat (but not necessarily faithfully flat) ring monomorphims, i.e., every flat ring monomorphism $\phi: R \into S$ induces a homomorphism $C^*_{\rho}(R \rtimes R\reg) \to C^*_{\rho}(S \rtimes S\reg)$ determined by $\rho(b,a) \ma \rho(\phi(b),\phi(a))$.
\eprop
\bproof
The first part about left semigroup C*-algebras is a consequence of the discussion in \cite[\S~2.5]{Li2}, where it was shown that the full semigroup C*-algebras $C^*(R \rtimes R\reg)$ from \cite[Definition~2.2]{Li2} are functorial for faithfully flat ring monomorphisms. At this point, we remark that a ring monomorphism $\phi: R \into S$ is faithfully flat if and only if $S / \phi(R)$ is flat over $R$ (see \cite[Chapitre~I, \S~3.5, Proposition~9]{Bour1}). As $R \rtimes R\reg \subseteq Q \rtimes Q\reg$ is left Toeplitz, because $R \rtimes R\reg$ satisfies the left independence condition by assumption and since $Q \rtimes Q\reg$ is amenable, \cite[Theorem~6.1]{Li3} tells us that full and reduced versions of our semigroup C*-algebras for $R \rtimes R\reg$ coincide, i.e., $C^*_s(R \rtimes R\reg) \cong C^*_{\lambda}(R \rtimes R\reg)$. As the two full versions $C^*(R \rtimes R\reg)$ and $C^*_s(R \rtimes R\reg)$ coincide by \cite[\S~3.1]{Li2}, we are done.

For the part about right semigroup C*-algebras, first of all observe that since $R \rtimes R\reg$ is right amenable, we have $C^*_s((R \rtimes R\reg)^{\text{op}}) \cong C^*_{\lambda}((R \rtimes R\reg)^{\text{op}}) \cong C^*_{\rho}(R \rtimes R\reg)$ by \cite[\S~4.1]{Li2}. So given a flat ring monomorphism $\phi: R \into S$, we can use the universal property of $C^*_s((R \rtimes R\reg)^{\text{op}})$ to construct the desired homomorphism. For $I \in \cI(R)$, let $\phi(I)S$ be the smallest ideal of $S$ generated by $\phi(I)$. Let $E_{S \times (\phi(I)S)\reg}$ be the orthogonal projection onto $\ell^2(S \times (\phi(I)S)\reg) \subseteq \ell^2(S \rtimes S\reg)$. We claim that the isometries $\rho(\phi(b),\phi(a))$, $(b,a) \in R \rtimes R\reg$, and the projections $E_{S \times (\phi(I)S)\reg}$, $I \in \cI(R)$, and $E_{\emptyset} = 0$, satisfy the relations I, II and III from \cite[Definition~3.2]{Li2} (the semigroup $P$ is $(R \rtimes R\reg)^{\text{op}}$ in our case). The first two relations are obviously fulfilled, so we only have to prove III. Let $Q(R)$ and $Q(S)$ be the quotient fields of $R$ and $S$. Assume that for $(b_i,a_i)$ and $(d_i,c_i)$ in $R \rtimes R\reg$, we have $(b_1,a_1)^{-1} (d_1,c_1) \dotsm (b_n,a_n)^{-1} (d_n,c_n) = (0,1)$ in $Q(R) \rtimes Q(R)\reg$. Then $(\phi(b_1),\phi(a_1))^{-1} (\phi(d_1),\phi(c_1)) \dotsm (\phi(b_n),\phi(a_n))^{-1} (\phi(d_n),\phi(c_n)) = (0,1)$ holds in $Q(S) \rtimes Q(S)\reg$. As relation III holds in $C^*_{\rho}(S \rtimes S\reg) \cong C^*_{\lambda}((S \rtimes S\reg)^{\text{op}})$, we have
\bglnoz
  && \rho(\phi(b_1),\phi(a_1))^* \rho(\phi(d_1),\phi(c_1)) \dotsm \rho(\phi(b_n),\phi(a_n))^* \rho(\phi(d_n),\phi(c_n)) \\
  &=& E_{(S \rtimes S\reg)(\phi(b_1),\phi(a_1)) (\phi(d_1),\phi(c_1))^{-1} \dotsm (\phi(b_n),\phi(a_n)) (\phi(d_n),\phi(c_n))^{-1}} \\
  &=& E_{S \times (\phi(c_n)^{-1} \phi(a_n) \dotsm \phi(c_1)^{-1} \phi(a_1) S)\reg}.
\eglnoz
Since we assume that $\phi$ is flat, we know that \cite[Lemma~2.18, (b)]{Li2} holds by \cite[Chapter~I, \S~2.6, Proposition~6]{Bour1}, so that equations (24) and (25) in \cite{Li2} hold as well. Therefore, $\phi(c_n)^{-1} \phi(a_n) \dotsm \phi(c_1)^{-1} \phi(a_1) S = \phi(c_n^{-1} a_n \dotsm c_1^{-1} a_1 R) S$, and thus,
\bglnoz
  && \rho(\phi(b_1),\phi(a_1))^* \rho(\phi(d_1),\phi(c_1)) \dotsm \rho(\phi(b_n),\phi(a_n))^* \rho(\phi(d_n),\phi(c_n)) \\
  &=& E_{S \times (\phi(c_n^{-1} a_n \dotsm c_1^{-1} a_1 R) S)\reg}.
\eglnoz
By universal property of $C^*_s((R \rtimes R\reg)^{\text{op}})$, there exists a homomorphism
\bglnoz
  && C^*_{\rho}(R \rtimes R\reg) \cong C^*_s((R \rtimes R\reg)^{\text{op}}) \to C^*_{\rho}(S \rtimes S\reg) \\ 
  && \rho(b,a) \ma \rho(\phi(b),\phi(a)), \, E_{R \times I\reg} \ma E_{S \times (\phi(I)S)\reg}.
\eglnoz
This means that we are done.
\eproof

In particular, since the inclusion $R \into Q$ is flat by \cite[Chapitre~II, \S~2.4, Th\'{e}or\`{e}me~1]{Bour1}, we have the following
\bcor
\label{rightC-->quot}
If $R$ satisfies the independence condition, then there exists a homomorphism $C^*_{\rho}(R \rtimes R\reg) \to C^*_{\rho}(Q \rtimes Q\reg) \cong C^*_{\lambda}(Q \rtimes Q\reg)$, $\rho(b,a) \ma \rho^{Q \rtimes Q\reg}_{(b,a)} \ma \lambda^{Q \rtimes Q\reg}_{(b,a)^{-1}}$.
\ecor

\section{K-theory}
\label{sec-K}

For integral domains which satisfy the independence condition, we now compute K-theory for $C^*_{\lambda}(R \rtimes R\reg)$ and $C^*_{\rho}(R \rtimes R\reg)$.
\btheo
\label{K}
Let $R$ be an integral domain with quotient field $Q$. Assume that $\cI(R)$ is independent. For $I \in \cI(R \subseteq Q)$, let $Q\reg_I \defeq \menge{a \in Q\reg}{aI = I}$ and $(R:I) = \menge{x \in Q}{xI \subseteq R}$. Moreover, $Q\reg$ acts on $\cI(R \subseteq Q)$ by (left and right) multiplication, and we let $Q\reg \backslash \cI(R \subseteq Q)$ and $\cI(R \subseteq Q)/Q\reg$ be the corresponding sets of orbits. Then
$$
  K_*(C^*_{\lambda}(R \rtimes R\reg)) \cong \bigoplus_{[I] \in Q\reg \backslash \cI(R \subseteq Q)} K_*(C^*_{\lambda}(I \rtimes Q\reg_I)),
$$
$$
  K_*(C^*_{\rho}(R \rtimes R\reg)) \cong \bigoplus_{[I] \in \cI(R \subseteq Q)/Q\reg} K_*(C^*_{\rho}((R:I) \rtimes Q\reg_I)).
$$
\etheo
Of course, it does not matter for commutative rings whether we write $Q\reg \backslash \cI(R \subseteq Q)$ or $\cI(R \subseteq Q)/Q\reg$, and the distinction between left and right group C*-algebras in our K-theoretic formulas is also not really necessary.
\bproof
We have already seen that $R \rtimes R\reg \subseteq Q \rtimes Q\reg$ is both left and right Toeplitz (see Lemma~\ref{Toeplitz}). Under our assumption that $\cI(R)$ is independent, we also know by Lemma~\ref{independence} that $R \rtimes R\reg$ satisfies both the left and right independence conditions. And finally, our group $Q \rtimes Q\reg$ is solvable, hence amenable, so that it satisfies the Baum-Connes conjecture with arbitrary coefficients by \cite{H-K}. Therefore, \cite[Corollary~4.9]{C-E-L2} applies in our situation (see also \cite[Remark~4.8]{C-E-L2}).

In order to deduce our K-theoretic formulas, all we have to do is to compute the orbits and stabilizer groups in our situation. For the left semigroup C*-algebras, two non-empty elements $(x+I) \times I\reg$ and $(y+J) \times J\reg$ of $\cJ_{\lambda}(R \rtimes R\reg \subseteq Q \rtimes Q\reg)$ lie in the same $Q \rtimes Q\reg$-orbit if and only if there exists $(b,a) \in Q \rtimes Q\reg$ such that $(b,a)((x+I) \times I\reg) = (y+J) \times J\reg$. The latter holds if and only if there exists $a \in Q\reg$ and $b \in Q$ with $aI=J$ and $b+ax+aI = y+J$. Since we can always find a suitable $b \in Q$, we see that $(x+I) \times I\reg$ and $(y+J) \times J\reg$ lie in the same $Q \rtimes Q\reg$-orbit if and only if there exists $a \in Q\reg$ and $b \in Q$ with $aI=J$. Thus, the map $Q \rtimes Q\reg \backslash \cJ_{\lambda}(R \rtimes R\reg \subseteq Q \rtimes Q\reg) \to Q\reg \backslash \cI(R \subseteq Q)$, $[(x+I) \times I\reg] \ma [I]$ is a bijection with inverse $[I \times I\reg] \mafr [I]$. Moreover, for $I \in \cI(R \subseteq Q)$, the left stabilizer group $\menge{(b,a) \in Q \rtimes Q\reg}{(b,a)(I \times I\reg) = I \times I\reg}$ is given by $I \rtimes Q\reg_I$ since
\bgloz
  (b,a)(I \times I\reg) = I \times I\reg \LRarr aI=I \ \wedge \ b+aI=I \LRarr b \in I \ \wedge \ a \in Q\reg_I.
\egloz
Plugging these observations into the formula from \cite[Corollary~4.9]{C-E-L2}, this proves our formula for $K_*(C^*_{\lambda}(R \rtimes R\reg))$. Analogous computations give the formula for $K_*(C^*_{\rho}(R \rtimes R\reg))$.
\eproof

\bcor
Let $R$ be a Krull ring with (multiplicative) units $R^*$ and divisor class group $C(R)$. For every $\mfk \in C(R)$, let $\mfa_\mfk$ be a divisorial ideal which represents $\mfk$. Then
$$
  K_*(C^*_{\lambda}(R \rtimes R\reg)) \cong \bigoplus_{\mfk \in C(R)} K_*(C^*(\mfa_\mfk \rtimes R^*)),
$$
$$
  K_*(C^*_{\rho}(R \rtimes R\reg)) \cong \bigoplus_{\mfk \in C(R)} K_*(C^*(\mfa_{\mfk^{-1}} \rtimes R^*)).
$$
\ecor
\bproof
For a Krull ring $R$, Lemma~\ref{Krull-independence} tells us that $\cI(R)$ is independent. So our theorem applies. As we have explained in \S~\ref{sec-Krull}, for a Krull ring, we have $\cI(R \subseteq Q) = D(R)$, so that $Q\reg \backslash \cI(R \subseteq Q) \cong C(R) \cong \cI(R \subseteq Q) / Q\reg$. Moreover, we know that the inverse in $(D(R),\bullet)$ of a divisorial ideal $\mfa$ is represented by $(R:\mfa)$, i.e., $(R:\mfa)$ represents the class $[\mfa]^{-1} \in C(R)$. And finally, given a divisorial ideal $\mfa$ and $a \in Q\reg$ with $a \mfa = \mfa$, we deduce that $(aR) \bullet \mfa = \mfa$ holds in $D(R)$. But since $D(R)$ is a group, this implies that $aR$ is the neutral element in $D(R)$, i.e., $aR=R$, and hence $a$ must lie in $R^*$. This shows that for every divisorial ideal $\mfa$ of $R$, we have $Q\reg_\mfa = R^*$. Plugging all these observations into the K-theoretic formulas from our previous theorem, we arrive at the desired result.
\eproof

\bremark
We see in particular that for a Krull ring, $K_*(C^*_{\lambda}(R \rtimes R\reg)) \cong K_*(C^*_{\rho}(R \rtimes R\reg))$. This is very surprising as in the following two sections, we will see that the left and right semigroup C*-algebras are completely different in general. In the context of Dedekind domains, this observation that the C*-algebraic difference between left and right semigroup C*-algebras becomes invisible in K-theory was already made in \cite[\S~6.4]{C-E-L2}.

At this point, a natural question would be: Does there exist an integral domain $R$ for which $K_*(C^*_{\lambda}(R \rtimes R\reg)) \ncong K_*(C^*_{\rho}(R \rtimes R\reg))$? Even for integral domains satisfying the independence condition, the answer to this question is not known. But looking at our general K-theoretic formulas, the only difference we can see is the difference between the ideal $I$ and the ideal $(R:I)$. So a more concrete question would be whether this difference in our formulas can lead to a different final outcome in K-theory, in the case of integral domains satisfying independence.
\eremark

\bremark
As we have seen in \S~\ref{sec-intdom}, another class of integral domains which satisfy the independence condition is given by integral domains which contain an infinite field. So also for these, the K-theoretic formulas from our theorem apply.
\eremark

\bremark
As in \cite[\S~8]{C-E-L1}, we can also compute K-theory for the semigroup C*-algebra of the multiplicative semigroup $R\reg$. As this is an abelian semigroup, we do not need to distinguish between left and right semigroup C*-algebras, and we simply write $C^*_r(R\reg) \defeq C^*_{\lambda}(R\reg) = C^*_{\rho}(R\reg)$. Let $R$ be a countable integral domain with quotient field $Q$ and group of multiplicative units $R^*$. If $R$ satisfies the independence condition, then, using the same notations as in Theorem~\ref{K}, we obtain the following K-theoretic formula:
$$
  K_*(C^*_r(R\reg)) \cong \bigoplus_{[I] \in Q\reg \backslash \cI(R \subseteq Q)} K_*(C^*(Q\reg_I)).
$$
In particular, if $R$ is a Krull ring, then we get
$$
  K_*(C^*_r(R\reg)) \cong \bigoplus_{\mfk \in C(R)} K_*(C^*(R^*)).
$$
\eremark

\bremark
As in \cite[\S~8]{C-E-L1}, all our K-theoretic formulas are actually formulas in KK-theory.
\eremark

\section{The primitive ideal space}
\label{sec-prim}

Let $R$ be a (countable) Krull ring with quotient field $Q$. First of all, let us describe the spectrum of $D \defeq D_{\lambda}(R \rtimes R\reg)$ (see \S~\ref{semiC}). Let $\cP \defeq \cP(R)$ be the set of prime ideals of $R$ which are of height $1$. For a function $m$: $\cP \to \Nz_0$, $\mfp \ma m_{\mfp}$ with finite support, we set $\cI_m \defeq \menge{\prod_{\mfp \in \cP} \mfp^{(v_{\mfp})}}{v_{\mfp} \leq m_{\mfp} \fa \mfp \in \cP}$. Recall that we take the product in $D(R)$, the monoid of divisorial ideals of $R$. Let $D_m \defeq C^*(\menge{E_{(r+I) \times I\reg}}{r \in R, \, I \in \cI_m})$. For $I \in \cI_m$ and $r+I \in R/I$, let $\chi_{r+I}$ be the character of $D_m$ determined by $\chi_{r+I}(E_{(s+J) \times J\reg}) = 1$ $\LRarr$ $r+I \subseteq s+J$ for $s \in R$, $J \in \cI_m$. Moreover, given $s \in R$, $J \in \cI_m$, set $S_{s+J} \defeq \menge{r+I \in \coprod_{I \in \cI_m} R/I}{r+I \subseteq s+J}$. Equip the set $\coprod_{I \in \cI_m} R/I$ with the topology determined by the following notion of convergence:

A sequence $(r_i+I_i)_i$ converges to $r+I$ in $\coprod_{I \in \cI_m} R/I$ if and only if for all $s+J \in \coprod_{I \in \cI_m} R/I$, there is $i_0 \in \Nz$ such that for all $i \geq i_0$,
$
r_i+I_i
\bfa
  \in S_{s+J} \falls r+I \in S_{s+J} \\
  \notin S_{s+J} \falls r+I \notin S_{s+J}.
\efa
$
A straightforward computation shows
\blemma
If we equip $\Spec(D_m)$ with the topology of pointwise convergence and $\coprod_{I \in \cI_m} R/I$ with the topology described above, then the map $\coprod_{I \in \cI_m} R/I \to \Spec(D_m)$, $r+I \ma \chi_{r+I}$ is a homeomorphism.
\elemma
Given two functions $m$, $n$: $\cP \to \Nz_0$, we say that $m \leq n$ if $m_{\mfp} \leq n_{\mfp}$ for all $\mfp \in \cP$. By construction, $D_m \subseteq D_n$ whenever $m \leq n$. Another direct computation yields the following
\blemma
Upon identifying $\Spec(D_m)$ with $\coprod_{I \in \cI_m} R/I$ and $\Spec(D_n)$ with $\coprod_{I \in \cI_n} R/I$, the map $\Spec(D_n) \to \Spec(D_m)$, $\chi \ma \chi \vert_{D_m}$ is given by
$$\pi_{n,m}: \, \coprod_{I \in \cI_n} R/I \to \coprod_{I \in \cI_m} R/I, \, r+I \ma r + \prod_{\mfp \in \cP} \mfp^{(\min(m_{\mfp},v_{\mfp}(I)))} \text{ (where } I = \prod_{\mfp \in \cP} \mfp^{(v_{\mfp}(I))} \text{)}.$$
\elemma
As $D = \ilim_m D_m$, we obtain
\bcor
$\Spec(D) \cong \plim_m \gekl{\coprod_{I \in \cI_m} R/I; \; \pi_{n,m}}$.
\ecor
In the following, let us write $\Omega$ for $\Spec(D)$. Now let us describe $\Omega_{\infty}$, the spectrum of $D_{\infty} \defeq D_{\lambda, R \rtimes R\reg \subseteq Q \rtimes Q\reg}$. Recall that $D_{\infty}$ is the smallest $Q \rtimes Q\reg$-invariant sub-C*-algebra of $\ell^{\infty} (Q \rtimes Q\reg)$ which contains $D$. In our particular situation, we have $D_{\infty} = C^*(\menge{E_{(q+J) \times J\reg}}{q \in Q, \, J \in \cI(R \subseteq Q)})$. Ordering $R\reg$ by divisibility, we have $D_{\infty} = \ilim_{a \in R\reg} E_{(a^{-1} R) \times a^{-1} R\reg} D_{\infty} E_{(a^{-1} R) \times a^{-1} R\reg}$. Let $\alpha$ be the action of $Q \rtimes Q\reg$ on $D_{\infty}$ given by $\alpha_{(b,a)} E_{(q+J) \times J\reg} = E_{(b+aq+aJ) \times aJ\reg}$. For $a \in R\reg$, we can identify $E_{(a^{-1} R) \times a^{-1} R\reg} D_{\infty} E_{(a^{-1} R) \times a^{-1} R\reg}$ with $D$ via $\alpha_{(0,a)}$ as $\alpha_{(0,a)} E_{(a^{-1} R) \times a^{-1} R\reg} = E_{R \rtimes R\reg}$ and $E_{R \rtimes R\reg} D_{\infty} E_{R \rtimes R\reg} = D$. We obtain continuous maps $\iota_a$: $\Omega \to \Omega_{\infty}$, where $\iota_a(\chi)(d) \defeq \chi \circ \alpha_{(0,a)} (E_{(a^{-1} R) \times a^{-1} R\reg} d E_{(a^{-1} R) \times a^{-1} R\reg}) = \chi(E_{R \rtimes R\reg} \alpha_{(0,a)}(d) E_{R \rtimes R\reg})$. Obviously, $\iota_a$ is injective, and we have $\Omega_{\infty} = \bigcup_{a \in R\reg} \iota_a(\Omega)$ as for every character $\chi$ of $D_{\infty}$, there exists $a \in R\reg$ with $\chi(E_{(a^{-1} R) \times a^{-1} R\reg}) = 1$. Finally, for $c \in R\reg$, let $\sigma_c$: $\Omega \to \Omega$ be given by $\sigma_c(\chi) \defeq \chi(\lambda(0,c)^* \sqcup \lambda(0,c))$. Obviously, $\sigma_c$ is injective, and we have for all $a$, $c$ in $R\reg$ that
\bgloz
\xymatrix{
\Omega \ar[d]_{\sigma_c} \ar[r]_{\iota_a} & \Omega_{\infty}
\\
\Omega \ar[ur]_{\iota_{ca}} &
}
\egloz
commutes:
\bglnoz
  \iota_{ca}(\sigma_c(\chi))(d) &=& \sigma_c(\chi)(E_{R \rtimes R\reg} \alpha_{(0,ca)}(d) E_{R \rtimes R\reg}) \\
  &=& \chi(\lambda(0,c)^* E_{R \rtimes R\reg} \alpha_{(0,a)}(d) E_{R \rtimes R\reg} \lambda(0,c)) \\
  &=& \chi(E_{R \rtimes R\reg} \alpha_{(0,c^{-1})}(\alpha_{(0,ca)}(d)) E_{R \rtimes R\reg}) \\
  &=& \chi(E_{R \rtimes R\reg} \alpha_{(0,a)}(d) E_{R \rtimes R\reg}) = \iota_a(\chi)(d)
\eglnoz
for all $\chi \in \Omega$ and $d \in D_{\infty}$. Putting all this together, we arrive at the following
\blemma
$\Omega_{\infty} \cong \ilim_{R\reg} \gekl{\Omega; \; \sigma_c}$ via the maps $\iota_a$: $\Omega \to \Omega_{\infty}$ ($a \in R\reg$) from above.
\elemma
\bcor
Let $\varphi$ be an element of $\iota_a(\Omega) \subseteq \Omega_{\infty}$ for $a \in R\reg$. For $C \subseteq \Omega_{\infty}$, we have $\varphi \in \overline{C}$ if and only if for all $m$: $\cP \to \Nz_0$ with finite support, there exists $\phi \in C \cap \iota_a(\Omega)$ with $\iota_a^{-1}(\varphi) \vert_{D_m} = \iota_a^{-1}(\phi) \vert_{D_m}$.
\ecor
Moreover, the right action $\alpha^*$ of $Q \rtimes Q\reg$ on $\Omega_{\infty}$ given by $\alpha^*_{(z,y)}(\chi) = \chi \circ \alpha_{(z,y)}$ can be described as follows: Assume that we have $(z,y) = (c^{-1}x,c^{-1}w)$ for $c \in R\reg$ and $(x,w) \in R \rtimes R\reg$. Then for $\chi \in \Omega$, $a \in R\reg$, we have $\alpha^*_{(z,y)}(\iota_a(\chi)) \in \iota_{aw}(\Omega)$, and if $\chi \vert_{D_m} = \chi_{r+I}$ holds for some function $m$: $\cP \to \Nz_0$ with finite support, we have $\iota_{aw}^{-1}(\alpha^*_{(z,y)}(\iota_a(\chi))) \vert_{D_{m+v(c)}} = \chi_{-ax+cr+cI}$, where $(m+v(c))_{\mfp} = m_{\mfp} + v_{\mfp}(c)$.

Our next goal is to describe the quasi-orbit space of $\alpha^*$. First of all, consider the following decomposition of $\cP$: We set $\cP_{\tfin} \defeq \menge{\mfp \in \cP}{[R:\mfp^{(i)}] < \infty \fa i \in \Nz}$ and $\cP_{\tinf} \defeq \cP \setminus \cP_{\tfin}$. We define $v_{\mfp}(\chi) \defeq \sup \menge{v_{\mfp}(I)}{\exists \, r \in Q \text{ with } \chi(E_{(r+I) \times I\reg}) = 1}$ for $\chi \in \Omega_{\infty}$ and $\mfp \in \cP$. Let $\cP(\chi) \defeq \menge{\mfp \in \cP}{v_{\mfp}(\chi) = \infty}$, $\cP_{\tfin}(\chi) \defeq \cP(\chi) \cap \cP_{\tfin}$ and $\cP_{\tinf}(\chi) = \cP(\chi) \cap \cP_{\tinf}$. Moreover, let $\chi \cdot (Q \rtimes Q\reg)$ be the orbit of $\chi$ under $\alpha^*$, i.e.,  $\chi \cdot (Q \rtimes Q\reg) = \menge{\chi \circ \alpha_{(b,a)}}{(b,a) \in Q \rtimes Q\reg}$.
\bprop
For $\chi \in \Omega_{\infty}$, we have
$$
  \overline{\chi \cdot (Q \rtimes Q\reg)} = \menge{\varphi \in \Omega_{\infty}}{\cP_{\tfin}(\chi) \subseteq \cP_{\tfin}(\varphi)}.
$$
\eprop
\bproof
Let us first show \an{$\supseteq$}. Take $\varphi \in \Omega_{\infty}$ with $\cP_{\tfin}(\chi) \subseteq \cP_{\tfin}(\varphi)$. Let $m$: $\cP \to \Nz_0$ be a map with finite support $\supp(m)$.

First of all, we claim that there exists $\psi \in \overline{\chi \cdot (Q \rtimes Q\reg)}$ with $\supp(m) \cap \cP(\psi) \subseteq \supp(m) \cap \cP(\varphi)$. To see this, take $\mfp \in \cP_{\tinf}(\chi)$. In the following, we construct $\psi \in \overline{\chi \cdot (Q \rtimes Q\reg)}$ with $\mfp \notin \cP(\psi)$ and $\cP(\psi) \subseteq \cP(\chi)$. This is all we need to show our claim because we have $\cP(\chi) \setminus \cP(\psi) \subseteq \cP_{\tinf}(\chi)$ so that an iteration of our construction yields $\psi \in \overline{\chi \cdot (Q \rtimes Q\reg)}$ with $\supp(m) \cap \cP(\psi) \subseteq \supp(m) \cap \cP(\varphi)$. Assume that $\chi$ lies in $\iota_a(\Omega)$ for $a \in R\reg$. Then $\chi \circ \alpha_{(0,a^{-1})}$ lies in $\Omega$, and we have $\cP(\chi \circ \alpha_{(0,a^{-1})}) = \cP(\chi)$. As $\mfp$ lies in $\cP(\chi)$, it lies in $\cP(\chi \circ \alpha_{(0,a^{-1})})$. As $\mfp$ lies in $\cP_{\tinf}$ and $\chi \circ \alpha_{(0,a^{-1})}$ lies in $\Omega$, there exists $i \in \Nz$, $r \in R$ with $[R:\mfp^{(i)}] = \infty$ and $\chi \circ \alpha_{(0,a^{-1})} (E_{(r+\mfp^{(i)}) \times {\mfp^{(i)}}\reg}) = 1$. Choose a sequence $(r_n)_n$ in $R$ such that $r-r_n+\mfp^{(i)}$ are pairwise disjoint. This is possible since $[R:\mfp^{(i)}] = \infty$. Then $\chi \circ \alpha_{(0,a^{-1})} \circ \alpha_{(r_n,1)} (E_{(r-r_n+\mfp^{(i)}) \times {\mfp^{(i)}}\reg}) = 1$. As the sequence $(\chi \circ \alpha_{(0,a^{-1})} \circ \alpha_{(r_n,1)})_n$ lies in $\Omega$ and since $\Omega$ is compact, there exists a convergent subsequence with limit $\psi \in \Omega$. By construction, $\psi$ lies in $\overline{\chi \cdot (Q \rtimes Q\reg)}$. Moreover, for every $s \in R$, $\psi(E_{(s+\mfp^{(i)}) \times {\mfp^{(i)}}\reg}) = 0$. Hence $\mfp \notin \cP(\psi)$. And for every $\mfq \in \cP$ and $v \in \Nz_0$, we know that $\chi \circ \alpha_{(0,a^{-1})}(E_{(r+\mfq^{(v)}) \times {\mfq^{(v)}}\reg}) = 0$ for all $r \in R$ implies that $\chi \circ \alpha_{(0,a^{-1})} \circ \alpha_{(r_n,1)}(E_{(r+\mfq^{(v)}) \times {\mfq^{(v)}}\reg}) = 0$ for all $r \in R$ and hence $\psi(E_{(r+\mfq^{(v)}) \times {\mfq^{(v)}}\reg}) = 0$ for all $r \in R$. This means that $v_{\mfq}(\psi) \leq v_{\mfq}(\chi \circ \alpha_{(0,a^{-1})})$ and thus $\cP(\psi) \subseteq \cP(\chi \circ \alpha_{(0,a^{-1})}) = \cP(\chi)$. Therefore, we have constructed $\psi \in \overline{\chi \cdot (Q \rtimes Q\reg)}$ with the desired properties.

Secondly, take $\psi \in \overline{\chi \cdot (Q \rtimes Q\reg)}$ with $\supp(m) \cap \cP(\psi) \subseteq \supp(m) \cap \cP(\varphi)$. Passing over to $\psi \circ \alpha_{(0,a^{-1})}$ if necessary, we may assume without loss of generality that $\psi$ lies in $\Omega$. For every $\mfp$ in $\supp(m)$ with $\mfp \notin \cP(\psi)$, there exists $r_{\mfp} \in R$ such that $\psi(E_{(r_{\mfp}+\mfp^{(v_{\mfp}(\psi))}) \times {\mfp^{(v_{\mfp}(\psi))}}\reg}) = 1$ and $\psi(E_{(r+\mfp^{(v)}) \times {\mfp^{(v)}}\reg}) = 0$ for all $v > v_{\mfp}(\psi)$ and $r \in R$. Let $\bigcap_{\substack{\mfp \in \supp(m) \\ \mfp \notin \cP(\psi)}} (r_{\mfp} + \mfp^{(v_{\mfp}(\psi))}) = s+J$ for some $s \in R$, $J \in \cI(R)$. Then $\psi(E_{(s+J) \times J\reg}) = 1$, and thus $\psi \circ \alpha_{(s,1)}(E_{J \times J\reg}) = 1$. Now assume that $\iota_a^{-1}(\varphi) \vert_{D_m} = \chi_{r+I}$. By Proposition~\ref{approx}, there exists $y$ in $Q\reg$ with
\bgln
\label{y:cond1}
  &v_{\mfp}(y)& = - v_{\mfp}(\psi) + v_{\mfp}(I) \fa \mfp \text{ in } \supp(m) \text{ with } \mfp \notin \cP(\psi), \\
\label{y:cond2}
  &v_{\mfp}(y)& \geq - v_{\mfp}(J) \fa \mfp \in \cP.
\egln
Then $\psi \circ \alpha_{(s,1)} \circ \alpha_{(0,y^{-1})} (E_{yJ \times yJ\reg}) = 1$. As $yJ \subseteq R$ by \eqref{y:cond2}, we know that $\psi \circ \alpha_{(s,1)} \circ \alpha_{(0,y^{-1})}$ lies in $\Omega$. We claim that $\psi \circ \alpha_{(s,1)} \circ \alpha_{(0,y^{-1})} \vert_{D_m} = \chi_{x+I}$ for some $x \in R$. To see this, write $\psi \circ \alpha_{(s,1)} \circ \alpha_{(0,y^{-1})} \vert_{D_m} = \chi_{x+I'}$ for some $x \in R$, $I' \in \cI_m$. For $\mfp$ in $\supp(m)$ with $\mfp \in \cP(\psi)$, we have $v_{\mfp}(I') = m_{\mfp} = v_{\mfp}(I)$. For $\mfp$ in $\supp(m)$ with $\mfp \notin \cP(\psi)$, we certainly have $v_{\mfp}(I') \geq v_{\mfp}(yJ) = v_{\mfp}(I)$, but we also have $v_{\mfp}(I') \leq v_{\mfp}(\psi \circ \alpha_{(s,1)} \circ \alpha_{(0,y^{-1})}) = v_{\mfp}(y) + v_{\mfp}(\psi \circ \alpha_{(s,1)}) = v_{\mfp}(y) + v_{\mfp}(\psi) = v_{\mfp}(I)$. Thus $I = I'$, i.e., $\psi \circ \alpha_{(s,1)} \circ \alpha_{(0,y^{-1})} \vert_{D_m} = \chi_{x+I}$. This implies that $\psi \circ \alpha_{(s,1)} \circ \alpha_{(0,y^{-1})} \circ \alpha_{(-r+x,1)} \vert_{D_m} = \chi_{x+I}$. Therefore, if we set
$$
  \phi \defeq \rukl{\psi \circ \alpha_{(s,1)} \circ \alpha_{(0,y^{-1})} \circ \alpha_{(-r+x,1)}} \circ \alpha_{(0,a)} 
  = \iota_a \rukl{\psi \circ \alpha_{(s,1)} \circ \alpha_{(0,y^{-1})} \circ \alpha_{(-r+x,1)}},
$$
then $\phi \in \iota_a(\Omega)$, $\phi \in \overline{\chi \cdot (Q \rtimes Q\reg)}$ and $\iota_a^{-1}(\phi) \vert_{D_m} = \psi \circ \alpha_{(s,1)} \circ \alpha_{(0,y^{-1})} \circ \alpha_{(-r+x,1)} \vert_{D_m} = \chi_{r+I} = \iota_a^{-1}(\varphi) \vert_{D_m}$. This shows that
$$
  \overline{\chi \cdot (Q \rtimes Q\reg)} \supseteq \menge{\varphi \in \Omega_{\infty}}{\cP_{\tfin}(\chi) \subseteq \cP_{\tfin}(\varphi)}.
$$
To prove \an{$\subseteq$}, assume that we have a sequence $\chi \circ \alpha_{(b_n,a_n)}$ in $\overline{\chi \cdot (Q \rtimes Q\reg)}$ with $\lim_{n \to \infty} \chi \circ \alpha_{(b_n,a_n)} = \varphi$. Without loss of generality, we can assume that $\varphi$ and $\chi \circ \alpha_{(b_n,a_n)}$, for all $n \in \Nz$, lie in $\Omega$. Every $\mfp \in \cP_{\tfin}(\chi)$ also lies in $\cP_{\tfin}(\chi \circ \alpha_{(b_n,a_n)}) = \cP_{\tfin}(\chi)$. Thus given $i \in \Nz$, there exists for every $n \in \Nz$ an element $r_n \in R$ with $\chi \circ \alpha_{(b_n,a_n)}(E_{(r_n + \mfp^{(i)}) \times {\mfp^{(i)}}\reg}) = 1$. But since $\mfp$ lies in $\cP_{\tfin}$, we must have $[R:\mfp{(i)}] < \infty$, so that there exists $r \in R$ with $\chi \circ \alpha_{(b_n,a_n)}(E_{(r + \mfp^{(i)}) \times {\mfp^{(i)}}\reg}) = 1$ for infinitely many $n \in \Nz$. Thus $\varphi(E_{(r + \mfp^{(i)}) \times {\mfp^{(i)}}\reg}) = 1$, and we conclude that $\mfp \in \cP_{\tfin}(\varphi)$.
\eproof
\bcor
Equip the power set $2^{\cP_{\tfin}}$ with the power-cofinite topology. Then we can identify the quasi-orbit space of $\alpha^*$ with $2^{\cP_{\tfin}}$ as topological spaces via $[\chi] \ma \cP_{\tfin}(\chi)$.
\ecor
Recall that the basic open sets in the power-cofinite topology on $2^{\cP_{\tfin}}$ are given by $\menge{S \subseteq \cP_{\tfin}}{F \cap S = \emptyset}$, $F \subseteq \cP_{\tfin}$ finite.

\bprop
\label{qo:trivstab}
Assume that $R$ is a (countable) Krull ring and that $\cP_{\tinf} \neq \emptyset$ or that $\cP_{\tfin}$ is infinite. Then for every subset $S \subseteq \cP_{\tfin}$, there exists $\chi \in \Omega_{\infty}$ with trivial stabilizer and $\cP_{\tfin}(\chi) = S$.
\eprop
We need a bit of preparation for the proof of this proposition:

By assumption, $R$ is countable, so that $Q$ and hence $Q \rtimes Q\reg$ are countable. Let $(b_1,a_1)$, $(b_2,a_2)$, ... be an enumeration of $Q \rtimes Q\reg \setminus \gekl{(0,1)}$, i.e., $Q \rtimes Q\reg \setminus \gekl{(0,1)} = \gekl{(b_i,a_i)}_{i \in \Nz}$. Also, $\cP$ is countable, so that $S$ is countable, and we let $\mfq_1$, $\mfq_2$, ... be an enumeration of $S$. This list of $\mfq_i$ may be finite. If $\cP_{\tinf} \neq \emptyset$, we just take an arbitrary $\mfp \in \cP_{\tinf}$ and set $\mfp_i \defeq \mfp$ for all $i \in \Nz$. If $\cP_{\tinf} = \emptyset$, then $\cP_{\tfin}$ must be infinite by assumption, and we let $\mfp_1$, $\mfp_2$, ... be an enumeration of $\cP_{\tfin}$. The key step in the proof of our proposition is the following
\blemma
Assume that we are in the situation of our proposition, and let the notation be as introduced above. Given $i \in \Nz$, $r \in Q$ and $I \in \cI(R \subseteq Q)$, there exists $s \in Q$ and $J \in \cI(R \subseteq Q)$ with
\begin{itemize}
\item $(s+J) \cap (b_i + a_i s + a_i J) = \emptyset$,
\item $s+J \subseteq r+I$,
\item $v_{\mfq_1}(J), \dotsc, v_{\mfq_i}(J) \geq i$ (if $i > \# S$, this just means that $v_{\mfq}(J) \geq i$ for all $\mfq \in S$),
\item $v_{\mfp}(J) = v_{\mfp}(I)$ for all $\mfp_i \neq \mfp \in \cP_{\tfin} \setminus S$
\end{itemize}
\elemma
\bproof
First, choose $J \in \cI(R \subseteq Q)$ with
\begin{itemize}
\item $J \subseteq I$,
\item $v_{\mfp_i}(J) > 
\bfa
  v_{\mfp_i}(b_i) \falls a_i = 1 \text{ (} \Rarr \, b_i \neq 0 \text{)}, \\
  v_{\mfp_i}(a_i - 1) + \abs{v_{\mfp_i}(a_i)} + v_{\mfp_i}(I) \falls a_i \neq 1,
\efa
$
\item $v_{\mfq_1}(J), \dotsc, v_{\mfq_i}(J) \geq i$,
\item $v_{\mfp}(J) = v_{\mfp}(I)$ for all $\mfp_i \neq \mfp \in \cP_{\tfin} \setminus S$.
\end{itemize}
If $a_i = 1$, $b_i \neq 0$, then since $v_{\mfp_i}(b_i) < v_{\mfp_i}(J)$, we must have $b_i \notin J$. In that case, we set $s \defeq r$, and we have $b_i + a_i s + a_i J = b_i + s + J \neq s + J$ and thus $(b_i + s + J) \cap (s + J) = \emptyset$. Moreover, $s+J \subseteq r+I$ obviously holds, and the remaining conditions in our lemma are also valid by construction.

Now assume that $a_i \neq 1$. Given $x$ and $y$ in $J$, we have
\bglnoz
  v_{\mfp_i}((a_i-1)^{-1}(x+a_i y)) &=& -v_{\mfp_i}(a_i - 1) + v_{\mfp_i}(x + a_i y) \\
  &\geq& -v_{\mfp_i}(a_i - 1) + \min(v_{\mfp_i}(x), v_{\mfp_i}(a_i) + v_{\mfp_i}(y)) \\
  &\geq& -v_{\mfp_i}(a_i - 1) - \abs{v_{\mfp_i}(a_i)} + v_{\mfp_i}(J) > v_{\mfp_i}(I).
\eglnoz
By Proposition~\ref{approx}, this shows that $I \nsubseteq (a_i - 1)^{-1}(J + a_i J)$. Thus $r+I \nsubseteq -(a_i - 1)^{-1} b_i + (a_i - 1)^{-1}(J + a_i J)$, so that we can choose $s \in r+J$ with $s \notin -(a_i - 1)^{-1} b_i + (a_i - 1)^{-1}(J + a_i J)$. The latter condition implies
\bglnoz
  && b_i + (a_i - 1)s \notin J + a_i J \\
  &\Rarr& b_i + a_i s + (J + a_i J) \neq s + (J + a_i J) \\
  &\Rarr& (b_i + a_i s + (J + a_i J)) \cap (s + (J + a_i J)) = \emptyset \\
  &\Rarr& (b_i + a_i s +  a_i J) \cap (s + J) = \emptyset.
\eglnoz
By construction, this choice of $s$ and $J$ satisfies all the desired conditions.
\eproof

\bproof[Proof of Proposition~\ref{qo:trivstab}]
Starting with $r_0 \defeq 0$, $I_0 \defeq R$, we can with the help of the previous lemma proceed inductively on $i$ to obtain a decreasing chain $R \supseteq r_1 + I_1 \supseteq r_2 + I_2 \supseteq \dotso$ such that
\begin{enumerate}
\item[I.] $(r_i + I_i) \cap (b_i + a_i r_i + a_i I_i) = \emptyset$ for all $i \in \Nz$,
\item[II.] $v_{\mfq_1}(I_i), \dotsc, v_{\mfq_i}(I_i) \geq i$ for all $i \in \Nz$,
\item[III.] $v_{\mfp_i}(I_j) = v_{\mfp_i(I_i)}$ for all $i, \, j \in \Nz$ with $j \geq i$ and $\mfp_i \in \cP_{\tfin} \setminus S$.
\end{enumerate}
Now let $\chi \in \Omega_{\infty}$ be the character uniquely determined by
$$\chi(E_{(s+J) \times J\reg}) = 1 \text{ if and only if there exists } i \in \Nz \text{ such that } r_i + I_i \subseteq s+J.$$
Existence of such a character is guaranteed by \cite[Corollary~2.8]{Li3} since our ring $R$ satisfies the independence condition by Lemma~\ref{Krull-independence}. By conditions II and III, we have $\cP_{\tfin}(\chi) = S$. Moreover, for all $i \in \Nz$, we must have $\chi \circ \alpha_{(b_i,a_i)} \neq \chi$. Otherwise, we would get $1 = \chi(E_{(r_i + I_i) \times I_i\reg}) =  \chi \circ \alpha_{(b_i,a_i)}(E_{(r_i + I_i) \times I_i\reg}) = \chi(E_{(b_i + a_i r_i + a_i I_i) \times a_i I_i\reg})$. But this contradicts condition I. This shows that $\chi$ has all the desired properties from our proposition.
\eproof

\bcor
\label{essfree}
In the situation of Proposition~\ref{qo:trivstab}, the action $\alpha^*$ of $Q \rtimes Q\reg$ on $\Omega_{\infty}$ is essentially free.
\ecor

\bcor
\label{primitive}
In the situation of Proposition~\ref{qo:trivstab}, there are order-preserving homeomorphisms $\Prim(C^*_{\lambda}(R \rtimes R\reg)) \cong \Prim(D_{\infty} \rtimes_{\alpha,r}(Q \rtimes Q\reg)) \cong 2^{\cP_{\tfin}}$, where $2^{\cP_{\tfin}}$ is given the power-cofinite topology. All the orders are given by inclusion (of primitive ideals or subsets).
\ecor
\bproof
For the first identification, use that $R \rtimes R\reg \subseteq Q \rtimes Q\reg$ is left Toeplitz and apply \cite[Corollary~3.10]{Li3}. The second identification follows from \cite[Corollary~2.7]{E-L} using the previous corollary.
\eproof

\bcor
In the situation of Proposition~\ref{qo:trivstab}, there exists a one-to-one correspondence between $\cP_{\tfin}$, the set of prime ideals $\mfp$ of $R$ with $\height(\mfp) = 1$ and $[R:\mfp^{(i)}] < \infty$ for all $i \in \Nz$, and the set of minimal non-zero primitive ideals of $C^*_{\lambda}(R \rtimes R\reg)$.
\ecor
\bremark
In the situation of Proposition~\ref{qo:trivstab}, it is possible to give a concrete description of the primitive ideals as in \cite[\S~3 and \S~4]{Li4}.
\eremark

\bremark
In the situation of Proposition~\ref{qo:trivstab}, assume that we have $\cP = \cP_{\tinf}$, i.e., for every $\mfp \in \cP$ there exists $i \in \Nz$ with $[R:\mfp^{(i)}] = \infty$. Then Corollary~\ref{primitive} tells us that $C^*_{\lambda}(R \rtimes R\reg)$ is simple. This is a first indication that the case of infinite quotients is somewhat special. Actually, a much stronger results holds, see Theorem~\ref{infquot}.
\eremark

\bremark
For the right semigroup C*-algebra $C^*_{\rho}(R \rtimes R\reg)$, the analogue of Corollary~\ref{essfree} cannot be true. The reason is that the canonical projection $C^*_{\rho}(R \rtimes R\reg) \to C^*_{\rho}(Q \rtimes Q\reg)$ (which exists by Corollary~\ref{rightC-->quot}) corresponds to a point in the spectrum of $D_{\rho, R \rtimes R\reg \subseteq Q \rtimes Q\reg}$ which is fixed by every element of $Q \rtimes Q\reg$.
\eremark

\section{Pure infiniteness}
\label{sec-pi}

Our goal in this section is to prove the following
\btheo
\label{pi}
Let $R$ be an integral domain. Suppose that $R$ is not a field, i.e., $R\reg \neq R^*$, and that the Jacobson radical $\cR$ of $R$ (the intersection of all maximal ideals of $R$) vanishes, i.e., $\cR=(0)$. Then $C^*_{\lambda}(R \rtimes R\reg)$ is purely infinite and satisfies the ideal property.
\etheo

With the help of \cite[Proposition~2.14]{P-R}, we obtain as an immediate consequence:
\bcor
Let $R$ be an integral domain which is not a field and has vanishing Jacobson radical. Then $C^*_{\lambda}(R \rtimes R\reg)$ is strongly purely infinite, i.e., $C^*_{\lambda}(R \rtimes R\reg) \cong C^*_{\lambda}(R \rtimes R\reg) \otimes \cO_{\infty}$.
\ecor

These C*-algebraic properties play an important role in the structure theory and the classification program for C*-algebras. The reader may find more information in \cite{K-R1}, \cite{K-R2}, \cite{P-R} and \cite{Ror}.

For the proof of this proposition, we first need a few preparations. We consider the general situation that we have a subsemigroup $P$ of a group $G$. As explained in \S~\ref{semiC}, the semigroup C*-algebra $C^*_{\lambda}(P)$ contains a canonical commutative sub-C*-algebra $D_{\lambda}(P)$. By \cite[Lemma~3.11]{Li2}, there exists a faithful conditional expectation $E: C^*_{\lambda}(P) \to D_{\lambda}(P)$ characterized by
$$E(\lambda(p_1)^* \lambda(q_1) \dotsm \lambda(p_n)^* \lambda(q_n)) 
= \bfa
  0 \falls p_1^{-1} q_1 \dotsm p_n^{-1} q_n \neq e \text{ in } G, \\
  \lambda(p_1)^* \lambda(q_1) \dotsm \lambda(p_n)^* \lambda(q_n) \sonst.
\efa
$$
Moreover, an ideal $\fI_D$ of $D_{\lambda}(P)$ is called invariant if for all $p \in P$, $\lambda(p) \fI_D \lambda(p)^* \subseteq \fI_D$.
\blemma
Let $\fI_D$ be an invariant ideal of $D_{\lambda}(P)$, and let $x \in C^*_{\lambda}(P)$ be a positive element. If $P \subseteq G$ is left Toeplitz and $G$ is exact, then $E(x) \in \fI_D$ implies that $x \in \spkl{\fI_D}$. Here $\spkl{\fI_D}$ is the ideal of $C^*_{\lambda}(P)$ generated by $\fI_D$.
\elemma
\bproof
Let $D_{\lambda, \ping}$ be as in \S~\ref{semiC}. As explained in \S~\ref{semiC}, we think of $C^*_{\lambda}(P)$ as a sub-C*-algebra of $D_{\lambda, \ping} \rtimes_r G$. Let $E^G: D_{\lambda, \ping} \rtimes_r G \to D_{\lambda, \ping}$ be the canonical conditional expectation. It is obvious that the conditional expectation $E$ from above is the restriction of $E^G$ to $C^*_{\lambda}(P)$. Moreover, let $\fI_{\ping}$ be the smallest $G$-invariant ideal of $D_{\lambda, \ping}$ which contains $\fI_D$. Let $E^G_{\fI}$ be the canonical faithful conditional expectation $(D_{\lambda, \ping}/\fI_{\ping}) \rtimes_r G \to D_{\lambda, \ping} / \fI_{\ping}$, and let $\pi: D_{\lambda, \ping} \to D_{\lambda, \ping} / \fI_{\ping}$ and $\pi^G: D_{\lambda, \ping} \rtimes_r G \to (D_{\lambda, \ping}/\fI_{\ping}) \rtimes_r G$ be the canonical projections. Then the following diagram commutes:
\bgl
\label{CD-condexp}
  \xymatrix@C=1mm{
  D_{\lambda, \ping} \rtimes_r G
  \ar[d]_{E^G} \ar[rr]^{\pi^G} 
  & \ \ \ \ \ \ \ \ \ \ \ \ & (D_{\lambda, \ping}/\fI_{\ping}) \rtimes_r G 
  \ar[d]^{E^G_{\fI}} \\
  D_{\lambda, \ping} 
  \ar[rr]^{\pi} 
  & \ \ \ \ \ \ \ \ \ \ \ \ & D_{\lambda, \ping}/\fI_{\ping}
  }
\egl
Now assume $x \in C^*_{\lambda}(P)_+$ satisfies $E(x) \in \fI_D$. Then $\pi(E^G(x)) = 0$. As \eqref{CD-condexp} commutes, we conclude that $E^G_{\fI}(\pi^G(x)) = 0$. But $E^G_{\fI}$ is faithful, and so we conclude that $\pi^G(x) = 0$. By assumption, $G$ is exact, so that $x \in \ker(\pi^G) = \fI_{\ping} \rtimes_r G$. On the whole, we obtain that $x$ lies in $E_P (\fI_{\ping} \rtimes_r G) E_P = \spkl{\fI_D}$. The last equality is \cite[Lemma~7.7]{Li3}, using the assumption that $\ping$ is left Toeplitz.
\eproof
\bcor
\label{ax+b-condexp}
Let $R$ be an integral domain, $\fI$ an ideal of $C^*_{\lambda}(R \rtimes R\reg)$ and $E$ the faithful conditional expectation $C^*_{\lambda}(R \rtimes R\reg) \to D_{\lambda}(R \rtimes R\reg)$ from above. Then, for a positive element $x$ in $C^*_{\lambda}(R \rtimes R\reg)$, $E(x) \in \fI$ implies that $x \in \fI$.
\ecor
\bproof
We just have to apply our previous lemma to $P = R \rtimes R\reg$, $G = Q \rtimes Q\reg$ ($Q$ is the quotient field of $R$) and $\fI_D = \fI \cap D_{\lambda}(R \rtimes R\reg)$. The conditions are satisfied since $R \rtimes R\reg \subseteq Q \rtimes Q\reg$ is left Toeplitz by Lemma~\ref{Toeplitz}, $Q \rtimes Q\reg$ is exact (even amenable) and $\fI \cap D_{\lambda}(R \rtimes R\reg)$ is obviously invariant. Thus if $x \in C^*_{\lambda}(R \rtimes R\reg)_+$ satisfies $E(x) \in \fI$, then $E(x)$ lies in $\fI \cap D_{\lambda}(R \rtimes R\reg) = \fI_D$, and the previous lemma tells us that $x$ lies in $\spkl{\fI_D} = \spkl{\fI \cap D_{\lambda}(R \rtimes R\reg)} \subseteq \fI$.
\eproof

This corollary leads to the following observation:
\bprop
\label{crit-pi}
Let $R$ be an integral domain. Suppose $R$ has the following property:

Given $I \in \cI(R)$ and $(b_i,a_i)$, $(b'_i,a'_i)$ in $R \rtimes R\reg$ for $1 \leq i \leq m$ with $(b_i,a_i) \neq (b'_i,a'_i)$ for all $1 \leq i \leq m$, we can always find $(b,a) \in R \rtimes R\reg$, $c \in R\reg$ and $r_1, \, r_2 \in R$ such that
\begin{enumerate}
\item[1${}_b$.] $b \in I$,
\item[2${}_b$.] $(b'_i-b_i)+(a'_i-a_i)b \neq 0$ for all $1 \leq i \leq m$,
\item[1${}_a$.] $a \in 1+I$,
\item[2${}_a$.] $(b'_i-b_i)+(a'_i-a_i)b \notin aR$ for all $1 \leq i \leq m$,
\item[($*_c$)] $c \notin R^*$ and $c \in 1+I$,
\item[($*_r$)] $r_1$, $r_2$ lie in $I$, and $r_1 + cR \neq r_2 + cR$.
\end{enumerate}

Then $C^*_{\lambda}(R \rtimes R\reg)$ is purely infinite and satisfies the ideal property.
\eprop
\bproof
First of all, Corollary~\ref{ax+b-condexp} tells us that the \cite[Corollary~8.2.8]{C-E-L1} holds for arbitrary integral domains. Secondly, going through the proofs of \cite[Lemma~8.2.9 and Lemma~8.2.10]{C-E-L1}, we see that the condition on our ring in our proposition is the only property of a Dedekind domain that we need in these proofs in \cite{C-E-L1}. Therefore, for every ring satisfying the hypothesis in our proposition, the analogues of \cite[Lemma~8.2.9 and Lemma~8.2.10]{C-E-L1} hold because their proofs carry over. And since the proof of \cite[Theorem~8.2.4]{C-E-L1} only uses \cite[Corollary~8.2.8, Lemma~8.2.9 and Lemma~8.2.10]{C-E-L1}, the same arguments as in \cite{C-E-L1} work for the rings which satisfy our hypothesis.
\eproof

\bproof[Proof of Theorem~\ref{pi}]
All we have to do is to verify the hypothesis of our proposition for an integral domain $R$ which is not a field and has zero Jacobson radical. Suppose we are given $I \in \cI(R)$ and $(b_i,a_i)$, $(b'_i,a'_i)$ in $R \rtimes R\reg$ for $1 \leq i \leq m$ with $(b_i,a_i) \neq (b'_i,a'_i)$ for all $1 \leq i \leq m$.

First of all, as $R\reg \neq R^*$, every ideal in $\cI(R)$ contains infinitely many elements. As 2${}_b$ only excludes finitely many possibilites for $b$, we can always find $b \in R$ satisfying 1${}_b$ and 2${}_b$. Now set $w \defeq \prod_{i=1}^m ((b'_i-b_i)+(a'_i-a_i)b)$. By our choice of $b$, we know that $w \neq 0$. We claim that $1 + (wR) \cap I \nsubseteq R^*$. Assume the contrary, i.e., $1 + (wR) \cap I \subseteq R^*$. Let $x$ be a non-zero element of $(wR) \cap I$. Then $1 + Rx \subseteq R^*$ implies by \cite[Proposition~1.9]{A-M} that $x$ lies in the Jacobson radical $\cR$ of $R$. This contradicts our assumption that $\cR = (0)$. Hence $1 + (wR) \cap I \nsubseteq R^*$. This also shows that $1 + (wR) \cap I \nsubseteq R^* \cup \gekl{0}$, because if $0$ lies in $1 + (wR) \cap I$, then $-1$ lies in $(wR) \cap I$, and thus $(wR) \cap I = R$, so that $1 + (wR) \cap I = R \nsubseteq R^* \cup \gekl{0}$ as $R\reg \neq R^*$ by assumption. Therefore, we can find an element $a \in R\reg \setminus R^*$ that lies in $1 + (wR) \cap I$. This element obviously satisfies 1${}_a$. In addition, we have $w \notin aR$ since otherwise, there would exist $y \in R$ with $w = ay$. At the same time, $a \in 1 + (wR) \cap I$ implies that $a = 1 + wz$ for some $z \in R$. Thus $a = 1 + ayz$ $\Rarr$ $a(1-yz) = 1$. But this contradicts $a \notin R^*$. So we must have $w \notin aR$. By construction of $w$, this ensures that 2${}_a$ holds as well.

It remains to find $c \in R\reg$ and $r_1, \, r_2 \in R$ satisfying ($*_c$) and ($*_r$). By the same argument as above, we know that $1+I \nsubseteq R^* \cup \gekl{0}$, so that we can choose $c \in R\reg \setminus R^*$ with $c \in 1+I$. Moreover, we have $I \nsubseteq cR$. Otherwise, we would have $c-1 \in I \subseteq cR$, i.e., there would exist $x$ in $R$ with $c-1 = cx$ $\Rarr$ $c(1-x) = 1$ $\Rarr$ $c \in R^*$ contradicting the choice of $c$. Hence $I \nsubseteq cR$. This implies that $(cR) \cap I \nsubseteq I$, and thus $[I:(cR) \cap I] > 1$. This means that we can choose two elements $r_1$ and $r_2$ in $I$ with $r_1 + (cR) \cap I \neq r_2 + (cR) \cap I$, i.e., $r_1 - r_2 \notin (cR) \cap I$. But then, we must have $r_1 - r_2 \notin cR$, and we are done.
\eproof

\bremark
In contrast to this, for every integral domain $R$, the right semigroup C*-algebra $C^*_{\rho}(R \rtimes R\reg)$ is not purely infinite as it projects onto $C^*_{\rho}(Q \rtimes Q\reg)$ by Corollary~\ref{rightC-->quot}, hence onto the non-zero commutative C*-algebra $C^*_{\rho}(Q\reg)$.
\eremark

Here is another situation where left semigroup C*-algebras of $ax+b$-semigroups are purely infinite.
\btheo
\label{infquot}
Let $R$ be a countable integral domain such that for every $I$ and $J$ in $\cI(R)$ with $I \subseteq J$, we have $[J:I] = \infty$. Then $C^*_{\lambda}(R \rtimes R\reg)$ is a UCT Kirchberg algebra. In particular, $C^*_{\lambda}(R \rtimes R\reg)$ is purely infinite and simple.
\etheo
\bproof
First of all, we have seen in \S~\ref{sec-intdom} that a ring $R$ as in our theorem satisfies the independence condition. As $Q \rtimes Q\reg$ is amenable ($Q$ is the quotient field of $R$), \cite[Theorem~6.1]{Li3} and \cite[\S~3.1]{Li2} together imply that $C^*_{\lambda}(R \rtimes R\reg) \cong C^*(R \rtimes R\reg)$, where $C^*(R \rtimes R\reg)$ is the full semigroup C*-algebra from \cite[Definition~2.2]{Li2}. But since we have $[J:I] = \infty$ for all $I$ and $J$ in $\cI(R)$ with $I \subseteq J$, it is straightforward to check that $C^*(R \rtimes R\reg)$ has the same universal property as the ring C*-algebra $\fA[R]$ from \cite[Definition~7]{Li1}. Now our claim follows from \cite[Corollary~3, Corollary~4 and Corollary~8]{Li1}.
\eproof
\bcor
If $R$ is a countable integral domain which contains an infinite field, then $C^*_{\lambda}(R \rtimes R\reg)$ is a UCT Kirchberg algebra.
\ecor
Again, these results have no analogues for the right semigroup C*-algebras of $ax+b$-semigroups.

\section{Integral domains which do not satisfy the independence condition}
\label{ind-fails}

We present an example of an integral domain for which the independence condition fails.
\bprop
Let $R$ be a noetherian integral domain. Assume that every non-zero factorial ideal of $R$ is divisorial, and that for every $a \in R\reg$, we have $[R:aR] < \infty$. Then $R$ satisfies the independence condition if and only if $R$ is integrally closed.
\eprop
\bproof
This follows immediately from \cite{B-Q}, Theorem~{2.6}.
\eproof
Note that a noetherian integral domain $R$ satisfies the property that every non-zero fractional ideal is divisorial if and only if $R$ is Gorenstein with Krull dimension at most $1$ (see \cite[Theorem~6.3]{Bass} and \cite[Corollary~4.3]{Hein}). Moreover, orders in number fields (or more generally, global fields) satisfy the finiteness condition for quotients, $[R:aR] < \infty$ for all $a \in R\reg$. Hence we obtain
\bcor
Let $R$ be an order in a number field (or global field). If every non-zero fractional ideal of $R$ is divisorial, then $R$ satisfies the independence condition if and only if $R$ is integrally closed.
\ecor
For instance, the order $R = \Zz[i \sqrt{3}]$ satisfies our condition that every non-zero fractional ideal is divisorial (see \cite[Example~4.2]{Ste}). Moreover, $R$ is not integrally closed. Therefore, $R$ does not satisfy independence. By Lemma~\ref{independence}, $R \rtimes R\reg$ does not satisfy neither the left nor the right independence condition. However, it is still possible to compute K-theory for the corresponding semigroup C*-algebras. This will be explained in \cite{Li-No}. Surprisingly, it again turns out that $K_*(C^*_{\lambda}(R \rtimes R\reg)) \cong K_*(C^*_{\rho}(R \rtimes R\reg))$ for our example $R = \Zz[i \sqrt{3}]$.

\end{document}